\documentclass[11pt]{article}
\usepackage{amssymb}

\advance\hoffset by -1.5cm
\advance\voffset by -1.5cm
\newtheorem{theorem}{Theorem}[section]

\newtheorem{e-proposition}[theorem]{Proposition}

\newtheorem{e-definition}[theorem]{Definition\rm}

\newcommand{\beq}{\begin{equation}}
\newcommand{\eeq}{\end{equation}}
\newcommand{\beqa}{\begin{eqnarray}}
\newcommand{\eeqa}{\end{eqnarray}}
\newcommand{\beqas}{\begin{eqnarray*}}
\newcommand{\eeqas}{\end{eqnarray*}}

\newcommand{\cK}{\mathcal{K}}

\newcommand{\cN}{\mathcal{N}}
\newcommand{\cP}{\mathcal{P}}

\newcommand{\real}{{\ensuremath{\mathbb{R}}}}
\newcommand{\dvg}{{\ensuremath{\rm div}}}
\newcommand{\pd}{\partial}
\newcommand{\sumnl}{\sum\nolimits}

\newcommand{\hsp}{\hspace{0.3mm}} 
\newcommand{\done}{\hfill $\Box$}

\newcommand{\myS}{{\rm S}}
\newcommand{\mys}{{\rm s}}
\newcommand{\myK}{{\rm K}}
\newcommand{\myk}{{\rm k}}
\newcommand{\HS}{{\rm H}}

\begin{document}

\begin{center}
\begin{LARGE}
{\bf Unbiased risk estimation and scoring rules}
\end{LARGE}
\vskip 4mm
\large{\sl Werner Ehm}
\end{center}

\begin{abstract}
Stein unbiased risk estimation is generalized twice, from the Gaussian shift model to nonparametric families of smooth densities, and from the quadratic risk to more general divergence type distances. The development relies on a connection with local proper scoring rules.
%\medskip
%\centerline{{\bf R\'{e}sum\'{e}}}
%La m\'{e}thode de l'estimation du risque Steinien est generalis\'{e}e doublement, du mod\`{e}le de translation Gaussien \`{a} des familles nonparam\'{e}triques, et du risque quadratique \`{a} des distances du type divergence plus g\'{e}n\'{e}raux. Le d\'{e}veloppement r\'{e}side dans une relation avec des r\`{e}gles d'estimation locales propres.
\end{abstract}

\section{Introduction: SURE and the Hyv\"arinen score}

Consider the problem of estimating the parameter $\theta$ in the standard Gaussian shift family $\ P_\theta = \cN(\theta,I_d),\ \theta \in \real^d,$ based on an observation $x \in \real^d.$ Let $T$ be an estimator of $\theta$ of the form $T = x + g(x).$ Using partial integration, Stein \cite{S} showed that under weak conditions about $g,$ the quadratic risk $R(T,\theta) = E_\theta\, |T - \theta|^2$ of $T$ can be estimated unbiasedly by the expression
$\widehat R(T) = 2\hsp \dvg\hsp g(x) + |g(x)|^2 + d\  $
called {\em SURE} (Stein unbiased risk estimate), so that $E_\theta\, \widehat R(T) = R(T,\theta)$ for every $\theta \in \real^d$. Here $| \cdot |$ and $\langle \cdot, \cdot \rangle$ denote the Euclidean norm and inner product on $\real^d,$ respectively, and $\dvg \hsp g$ is the divergence of $g.$   
If in particular $g = \nabla \log f$ for some function $f > 0$ on $\real^d$, the risk estimate becomes
\beq\label{Sure2}
\widehat R(T) = 2\hsp \Delta \log f(x) + |\nabla \log f(x)|^2 + d,
\eeq
where as usual, $\nabla$ denotes the gradient and $\Delta = \dvg \hsp \nabla$ the Laplace operator on $\real^d.$ 
This special case occurs if $T$ is the posterior mean with respect to a prior distribution $\pi$: then $T = x + \nabla \log f(x)$ where $f(x) = \int p_\theta(x)\, d\pi(\theta)$ is the corresponding mixture density, so that $g = \nabla \log f.$

The striking similarity between SURE and the {\em Hyv\"arinen score} \cite{H05}, 
\beq\label{HS}
\HS( \hsp p,x) = 2 \hsp\hsp \frac{\Delta \hsp p(x)}{p(x)} - \left| \frac{\nabla p(x)}{p(x)} \right|^2 =\,  2\hsp \Delta \log p(x) + |\nabla \log p(x)|^2,
\eeq
has been noted in, e.g., \cite{H08}. In eq.~(\ref{HS}), $p$ denotes a sufficiently smooth, strictly positive probability density on $\real^d$. Originally, the Hyv\"arinen score was introduced for score matching, a minimum distance type estimation method. 
Its formal similarity to SURE is substantiated on reexpressing the risk of $T$ as a distance between densities. Consider the {\em Hyv\"arinen divergence} defined for smooth, positive densitites $p,q$ on $\real^d$ as
\beq\label{Hloss1}
d_\HS(p,q)  = \int |\nabla \log p(y) - \nabla \log q(y)|^2\, q(y)\, dy.
\eeq
If $p=f$ is a mixture density as above and $q=p_\theta$ is the density of $P_\theta$, we have 
$
\nabla \log f(x) - \nabla \log p_\theta(x) = \nabla \log f(x)+ x -\theta = T - \theta,
$
where again $T = x +\nabla \log f(x)$ is the corresponding posterior mean. Consequently, 
\beq\label{parloss}
R(T,\theta) = E_\theta\, |T - \theta|^2 = \int |\nabla \log f(x) - \nabla \log p_\theta(x)|^2\, p_\theta(x)\, dx  = d_\HS(f, p_\theta),
\eeq
that is, {\em the risk $R(T,\theta)$ of the parameter estimate $T= x + \nabla \log f(x)$ equals a distance between densities,} $d_\HS(f, p_\theta).$ Furthermore, the analogue of SURE in the density scenario is the Hyv\"arinen score $\HS(f,x),$ essentially.
In fact, Hyv\"arinen's idea, reinventing Stein's, was to apply partial integration to (\ref{Hloss1}) which, assuming boundary terms vanish, gives 
\beq\label{Hloss2}
d_\HS(p,q) = \int \left( 2\Delta \log p(y) + |\nabla \log p(y)|^2 \right) q(y)\, dy + \int |\nabla \log q(y)|^2\, q(y)\, dy;
\eeq
cf.~\cite{DL}, \cite{H05}. Since $\int |\nabla \log p_\theta(x)|^2\, p_\theta(x)\, dx = d \ ( \theta \in \real^d)$ in the standard normal case, where $q = p_\theta,$ it follows that 
\beq\label{densloss}
E_\theta \left(\HS(f,x) +  d \right)  = \, E_\theta \left(2\Delta \log f(x) + |\nabla \log f(x)|^2 + d \right)  =  d_\HS(f, p_\theta).
\eeq
That is, {\em the modified Hyv\"arinen score $H(f,x)+d$ respresents an unbiased estimate of the distance $d_\HS(f, p_\theta)$ of $f$ from the unknown ``true" density $p_\theta$, for any density $f>0$ on $\real^d$ satisfying suitable regularity conditions.}

The purpose of this note is to expand on this aspect of unbiased risk estimation by tying it to scoring rules. Local proper scoring rules are constructed as gradients of concave functionals \cite{GR}, \cite{HB}, and then shown to generalize SURE in that they furnish unbiased estimates of modified Bregman type distances. The development is related to (parts of) work by Dawid and Lauritzen \cite{DL}. See also \cite{EG}, \cite{PDL}.

\section{Local proper scoring rules and unbiased risk estimation}

We restrict the discussion of scoring rules to the setting relevant for this note, and refer to \cite{GR} for general information. Let $\cP$ denote the class of all probabilitiy densities with respect to the Lebesgue measure on $\real^d$ such that the following conditions hold for every $p \in \cP$: 
(P1) $p \in C^2;$ $\! $
(P2) $p>0$ everywhere on $\real^d;\, $ 
(P3) for every $m>0$ and $i,j \in \{1,\ldots,d\}$
\[
\lim_{|x| \to \infty} \, |x|^m \left( p(x) + |\pd_{x_i} p(x)| + |\pd_{x_i x_j}^2 p(x)| \right) = 0;
\]
(P4) there exists $a = a(p)>0$ such that for $i,j \in \{1,\ldots,d\}$,
\[
\lim_{|x| \to \infty} \, |x|^{-a} \left(|\log p(x)| +  \left[\frac{\pd_{x_i} p(x)}{p(x)}\right]^2 + \frac{|\pd_{x_i x_j}^2 p(x)| }{p(x)}\right) = 0.
\]
The class $\cP$ is quite large, being convex and comprising, e.g., all normal and logistic distributions. 

A {\em scoring rule} is a mapping $\myS: \cP \times \real^d \to \real$ assigning a numerical score, $\myS(p,x),$ to the density forecast, $p,$ when the observation that materializes is $x.$ We write 
$\myS(p,q) = \int \myS(p,x) \, q(x)\, dx = E_q\, \myS(p,\cdot)$
for the expected score when the density forecast is $p$ and the probability measure underlying $x$ is $q(x) dx.$ The scoring rule $\myS$ is {\em (strictly) proper} relative to $\cP$ if $\myS(q,q) \le \myS(p,q)$ for all $p,q \in \cP$ (with equality only if $p=q$).
The scoring rule $\myS$ is {\em local} (of order two, for the class $\cP$) if there exists a real function $\mys$ such that 
\[
\myS(p,x) = \mys\!\left(x,\log p(x), \nabla \log p(x),\nabla^2 \log p(x)\right) \qquad (p \in \cP, \ x \in \real^d),
\]
$\nabla^2 f(x)$ denoting the Hessian matrix of second-order partial derivatives of a function $f: \real^d \to \real$ at $x.$

The classical example of a (strictly) proper local scoring rule is the logarithmic score, $\myS(p,x) = -\log p(x).$ 
Another example is the Hyv\"arinen score (\ref{HS}). The latter can be regarded as being local of order two, in the obvious sense, and the former as local of order zero. Local scoring rules of any order $m \ge 0$ were recently investigated in \cite{PDL}, in the case $d=1.$ Hereafter, ``local" always is understood as ``local of order two." 

The following result lifts the construction of local proper scoring rules in  \cite{EG} from the one- to the higher-dimensional case $d\ge 1.$ Let $\cK$ denote the class of the {\em kernels} $\myk: \real^d \times \real^d \to \real$ satisfying the following conditions: 
$ $ (K1) $\ \myk \in C^2;$ $\, $
(K2) $\ $there are constants $C,\, r \in (0,\infty)$ such that whenever $\myk^\ast$ stands for the function $\myk = \myk(x,y)$ or any of its partial derivatives up to order two, then 
$|\myk^\ast(x,y)| \ \le \ C  \left(1+|x| + |y| \right)^r \, (x, y \in \real^d).$ 
With any $\myk \in \cK$ we associate a functional $\Phi = \Phi_\myk : \cP \to \real$ defined by
\beq\label{Phidef}
\Phi(p) \, = \,\int\nolimits_{\real^d} \myk\!\left(x,\nabla \log p(x)\right)\, p(x)\, dx \qquad (p \in \cP). 
\eeq
In view of the growth and decay conditions (K2), (P4), and (P3), the integral in (\ref{Phidef}) exists and is finite for every $p \in \cP.$ 
Let $\nabla_{\! y} \hsp \myk$ denote the partial gradient referring to the argument $y \in \real^d$ of $\myk=\myk(x,y),$  and recall that $\dvg\hsp g(x)$ stands for the trace of the total derivative at $x$ of a function $x \mapsto g(x)$ mapping $\real^d$ into itself.

\begin{theorem} \label{thm1}
Let $\myk \in \cK,$ and suppose that the associated functional $\Phi$ is concave on $\cP.$ Then
\beq\label{srform}
\myS(p,x) \, = \, \myk(x,\nabla \log p(x))\, - \, \frac{1}{p(x)}\ \dvg \bigg[\, p(x)\,  \nabla_{\! y} \hsp \myk\left(x,\nabla \log p(x)\right)\! \bigg] 
\eeq
is a local proper scoring rule relative to $\cP.$  It is strictly proper if $\Phi$ is strictly concave. Furthermore, if $y \mapsto \myk(x,y)$ is concave on $\real^d$ for every $x \in \real^d,$ then the functional $\Phi$ is concave on $\cP.$
\end{theorem}

\noindent
The {\em Proof } follows similar lines as in the case $d=1,$ see \cite[Sections 4.1, 5.1]{EG}. We only indicate that the tangent construction in \cite[Section 4.1]{EG} yields the scoring rule (\ref{srform}). To compute the (weak) gradient of $\Phi$ at $q \in \cP,$ let $p_t = (1-t)q + t p$ where $p \in \cP,\ t \in [0,1].$ 
Formal differentation ignoring all technicalities gives
\begin{equation}  \label{gcalc1} 
\frac{d}{dt} \left[ \Phi(\hsp p_t) \right]  
= \int \frac{\pd}{\pd t} \left[ \myK_{p_t} \hsp p_t \right] dx = \int \left[\myK_{p_t}\right] \hsp (p-q)\, dx + \int \left[ \frac{\pd}{\pd t} \myK_{p_t} \right]\! p_t\, dx,
\end{equation} 
wherein we put $\myK_{p_t}(x) = \myk(x , \nabla \log p_t(x))$ and omitted the argument $x$ of the integrands. 
For the last integral in (\ref{gcalc1}) we get by the divergence theorem, assuming the boundary integral vanishes,
\beq \label{gcalc2}
\int \left\langle \nabla_{\! y} \hsp \myk\left(\cdot\, ,\nabla \log p_t\right), \nabla \left( \frac{p-q}{p_t} \right)\right\rangle p_t\, dx = -\int \dvg\hsp  \bigg[ p_t\,  \nabla_{\! y} \hsp \myk\left(\cdot\, ,\nabla \log p_t\right) \bigg]\, \frac{p-q}{p_t}\ dx.
\eeq
Setting $t = 0$ in (\ref{gcalc1}) and (\ref{gcalc2}) and noting that $p_0 = q$ we find that
\begin{equation}  \label{dirderiv} 
\frac{d}{dt} \left[ \Phi(\hsp p_t) \right] \bigg|_{t=0} 
= \int \left\{\myk(\cdot\, , \nabla \log q) - \frac{1}{q}\, \dvg\hsp  \bigg[ q\,  \nabla_{\! y} \hsp \myk\left(\cdot\, ,\nabla \log q\right)\! \bigg] \right\} (\hsp p - q)\, dx \, .
\end{equation}
Thus, the gradient of $\Phi$ at $q$ is given by the expression in curly brackets in (\ref{dirderiv}). The scoring rule resulting from the tangent construction, $\myS(q,\cdot),$ differs from this gradient only by a correction term which can be shown to vanish. The negligibility of the boundary integral in (\ref{gcalc2}), and all the technicalities (existence of integrals, exchangeability of differentiation and integration, etc.) can be settled similarly as in \cite[Section 4.1]{EG}, using the assumptions made about the classes $\cP$ and $\cK.$ \done

\smallskip 
Any convex combination of a scoring rule $\myS$ as in Theorem \ref{thm1} with the logarithmic score yields a local proper scoring rule. In the case $d=1,$ scoring rules of this form exhaust the class of all local proper scoring rules \cite{EG}, \cite{PDL}. The complete characterization in the case $d>1$ remains open. 

\smallskip\noindent
{\bf Examples.} Let $\myk \in \cK$ be a kernel of the form $\myk(x,y) = \myk(y) = \psi(|y|),$ where $\psi$ is a concave $C^2$-function on $[0, \infty)$ with $\psi(0)=\psi'(0)=0.$ Then $y \mapsto \myk(y)$ is concave on $\real^d,$ and the corresponding scoring rule (\ref{srform}) is proper. Explicitly we have
\[
\myS(p,\cdot\, ) = \psi(|\sigma|) - \frac{\psi'(|\sigma|)}{|\sigma|} \left(|\sigma|^2 + \Delta \log p\, \right) - \left[ \psi''(|\sigma|) - \frac{\psi'(|\sigma|)}{|\sigma|} \right] \left\langle \frac{\sigma}{|\sigma|}\, , \left(\nabla^2 \log p\, \right) \frac{\sigma}{|\sigma|} \right\rangle 
\]
where $\sigma = \nabla \log p.$ For $\psi(t) =  - t^2$ we obtain the Hyv\"arinen score (\ref{HS}); putting $\psi(t) = - \log \cosh t$ yields another interesting example parallel to \cite[Example 5.3]{EG}.

A local scoring rule $\myS$ that is proper relative to $\cP$ gives rise to a Bregman type divergence measure $d_\myS(p,q) = \myS(p,q) - \myS(q,q)$ on $\cP \times \cP.$ The following representation of $d_\myS$ is closely related to \cite[Eq.~(53)]{PDL}.

\begin{theorem} \label{thm2}
Suppose that $\myS$ is of the form {\em (\ref{srform})} for some kernel $\myk \in \cK$ such that $y \mapsto \myk(x,y)$ is concave on $\real^d$ for every $x \in \real^d.$ Then the divergence $d_\myS$ admits the representation
\beqa \label{d-explicit}
&& \!\!  d_\myS(p,q) \\ \!\! & = & \!\!  E_q \left\{\myk(\cdot\, , \nabla \log p) \, -\ \myk(\cdot\, , \nabla \log q) \, + \, \bigg\langle  \frac{\nabla q}{q} - \frac{\nabla p}{p},\,  \nabla_{\! y} \hsp \myk\left(\cdot,\nabla \log p\right) \bigg\rangle \right\} \quad (p,q \in \cP).\nonumber
\eeqa
\end{theorem}

\noindent
{\em Proof. } Let $p,\, q \in \cP.$ By the assumptions on $\cP$ and $\cK,$ the divergence theorem applied to the scalar function $u(x) = q(x)/p(x)$ and the vector field $v(x) = p(x)\, \nabla_{\! y} \hsp \myk\left(x,\nabla \log p(x)\right)$ gives
\beqa\label{partint}
&& \!\! \lim_{r \to \infty}\,  -\!\int_{|x| \le r}\frac{q}{p}\ \dvg\hsp  \bigg[ p\  \nabla_{\! y} \hsp \myk\left(\cdot,\nabla \log p\right)\! \bigg] dx \\ \!\! & = & \!\!  \lim_{r \to \infty} \int_{|x| \le r} \bigg\langle  \frac{\nabla q}{q} - \frac{\nabla p}{p},\,  \nabla_{\! y} \hsp \myk\left(\cdot,\nabla \log p\right) \bigg\rangle\, q\, dx . \nonumber
\eeqa
The relation (\ref{d-explicit}) follows on writing $d_\myS(p,q) = E_q \{\myS(p,\cdot\, ) - \myS(q,\cdot\, )\},$ substituting (\ref{srform}) and using (\ref{partint}), and observing that $\, \int q^{-1} \dvg \left( q\,  \nabla_{\! y} \hsp \myk\left(\cdot,\nabla \log q\right) \right) q\, dx = 0.$ \done

\smallskip
Note that the expression in curly brackets in (\ref{d-explicit}) is nonnegative because for a concave function $f$ on $\real^d$ one has $f(y_1) - f(y_2) \ge \langle y_1-y_2, \nabla f(y_1) \rangle\ (y_1,\, y_2 \in \real^d).$  For the Hyv\"arinen score, where $\myk(x,y) = -|y|^2,$ that expression becomes $|\nabla \log p - \nabla \log q|^2,$ and $d_\myS$ becomes the Hyv\"arinen divergence (\ref{Hloss1}). 

To clarify the connection with SURE we note that the partial integration in (\ref{partint}) was used conversely by Stein and Hyv\"arinen, to pass from the risk representation (\ref{d-explicit}) to an expression of the form $E_q \{\myS(p,\cdot\, ) - \myS(q,\cdot\, )\}.$ In the latter, the scoring rule $\myS(p,\cdot\, )$ may serve as an unbiased estimate of $E_q\, \myS(p,\cdot\, ),$ while the term $E_q\, \myS(q,\cdot\, )$ is the same for all candidates $p,$ hence can be ignored if the focus is on risk {\em comparison.}
In nonparametric density estimation, e.g., risk comparison of competing estimates is applied for bandwidth selection, using cross-validation. 
Briefly, if $\widehat p_n = \widehat p_n(\cdot\, |x_1,\ldots,x_n)$ is an estimate of the unkown density $q \in \cP$ underlying the i.~i.~d.~observations $x_1,\ldots,x_n$ that is symmetric in the $x_i$, the cross-validated expression 
$\widehat R_n(\widehat p_{n-1})  =  n^{-1}\, \sumnl_{i=1}^n \myS\!\left( \widehat p_{n,-i}\, ,\hsp x_i \right)$ 
is an unbiased estimate of $R_{n-1}(\widehat p_{n-1},q),$ where $R_n(\widehat p_n,q) = E_q\, \myS(\widehat p_n, q)$ denotes the modified risk ignoring the term $E_q\, \myS(q, \cdot\, ) = \myS(q,q),$ which depends only on $q.$

The possibility of risk estimation is of course not confined to the local scoring rules considered here, as any proper scoring rule $\myS,$ whether local or not, gives rise to a divergence measure $d_\myS.$  Therefore, cross-validatory estimation of the (modified) risk generally is feasible, although exact unbiasedness as with the local scoring rules may not be achievable when global terms are involved. For example, unbiased estimation of the term $\int p(x)^2\, dx$ entering the quadratic score \cite{GR} does not seem possible. 

The particular interest of the scoring rules of the form (\ref{srform}) ensues from the fact that they do not require the knowledge of the normalizing constants of the probability densities, which may be unknown or hard to obtain in complex settings \cite{H05}, \cite{PDL}. This advantage can be combined with other desirable features such as improved robustness by working, for instance, with the log cosh scoring rule mentioned above.

\section*{Acknowledgement}
The author thanks Tilmann Gneiting, Steffen Lauritzen, and Matthew Parry for comments on earlier drafts of the paper. 
%The author thanks Tilmann Gneiting and Matt Parry for helpful comments on earlier drafts of the paper. 

\medskip
\noindent
Author address:\\
Institute for Frontier Areas of Psychology and Mental Health\\
Wilhelmstr.~3a, 79098 Freiburg, Germany

\smallskip
\noindent
e--mail:  {\tt ehm@igpp.de}


\begin{thebibliography}{00}
% please try to use the bibitem system -
% the references should be in alphabetical order of authors' names.
% Articles with a single author first, author will 1 co-author next,
% then author with several co-authors;


% \bibitem{label}
% Text of bibliographic item



\bibitem{DL} A P Dawid and S L Lauritzen (2005). The geometry of decision theory. In {\em Proc.~2nd Int.~Symp.~Inf.~Geom.~Appl.} 22-28. Univ.~Tokyo

\bibitem{EG} W Ehm and T Gneiting  (2011).  Local proper scoring rules of order two.  Preprint, \begin{tt} {arXiv:1102.5031v1} \end{tt}

\bibitem{GR} T Gneiting and A E Raftery (2007). Strictly proper scoring rules, prediction, and estimation. {\em J.~Amer.~Statist.~Assoc.} {\bf 102}, 359-378

\bibitem{HB} A D Hendrickson and R J Buehler (1971).  Proper scores for probability forecasters.  {\em Ann.~Math.~Statist.} {\bf 42}, 1916-1921

\bibitem{H05} A Hyv\"arinen  (2005). Estimation of non-normalized statistical models using score matching. {\em J.~Mach.~Learn.~Res.} {\bf 6}, 695-709

\bibitem{H08} A Hyv\"arinen  (2008). Optimal approximation of signal priors. {\em Neural Computation} {\bf 20}, 3087-3110

\bibitem{PDL} M Parry, A P Dawid and S Lauritzen  (2011).  Proper local
scoring rules.  Preprint, \begin{tt} {arXiv:1101.5011v1} \end{tt}

\bibitem{S} C M Stein (1981). Estimation of the mean of a multivariate normal distribution. {\em Ann.~Statist.} {\bf 9}, 1135-1151

\end{thebibliography}
\end{document}